\pdfoutput=1

\documentclass[a4paper, 12pt]{article}
\usepackage{amsmath}
\usepackage{amsthm}
\usepackage{amssymb}
\usepackage[T1]{fontenc}
\usepackage{cancel}
\usepackage[sc]{mathpazo}
\usepackage{enumerate}
\usepackage{etoolbox}
\usepackage{xcolor}

\patchcmd{\section}{\scshape}{\bfseries}{}{}

\newcommand{\R}{{\mathbb{R}}}
\newcommand{\Q}{{\mathbb{Q}}}

\newcommand{\Z}{{\mathbb{Z}}}

\newcommand{\stateSpace}{\mathcal{S}}

\DeclareMathOperator{\tr}{tr}

\newtheorem{theorem}{Theorem}
\newtheorem{proposition}[theorem]{Proposition}
\newtheorem{lemma}[theorem]{Lemma}
\newtheorem{corollary}[theorem]{Corollary}

\theoremstyle{definition}

\theoremstyle{remark}

\newtheorem{example}[theorem]{Example}

\newcommand{\intersect}{\cap}

\newcommand{\Sym}{\mathrm{Sym}}

\newcommand{\tensor}{\otimes}

\newcommand{\dual}{\vee}

\renewcommand{\epsilon}{\varepsilon}

\newcommand{\X}{\mathcal{X}}

\newcommand{\innerProduct}[2]{\langle #1, #2 \rangle}

\newcommand{\K}{\mathbb{K}}

\newcommand{\isom}{\cong}

\DeclareMathOperator{\Psd}{\mathrm{Psd}}

\renewcommand{\matrix}[1]{\mathbf{#1}}
\renewcommand{\vector}[1]{\mathbf{#1}}

\renewcommand{\subset}{\subseteq}

\newcommand{\algebraicInterior}{\circ}

\newcommand{\Du}{Du'}

\newcommand{\Le}{L\^{e}}

\newcommand{\matrixSize}{q}
\newcommand{\numberOfConstraints}{m}

\newcommand{\degree}{e}
\newcommand{\constantBound}{n}

\newcommand{\nonnegativeLocus}{\mathcal{X}}
\newcommand{\numberOfModules}{m}
\newcommand{\dimension}{n}

\begin{document}

\providecommand{\keywords}[1]
{
  \small	
  \textit{Keywords:   } #1
}

\providecommand{\msccode}[1]
{
  \small	
  \textit{2020 MSC:   } #1
}

\title{Archimedean Representation Theorem for modules over a commutative ring}
\author
{Colin Tan\footnote{Email: \texttt{colinwytan@protonmail.com}} 
}

\date{06 Nov 2023}

\maketitle

\noindent \keywords{archimedean semiring, module over a commutative ring, 
order unit, polynomial, polytope, positive definite matrix, Positivstellensatz, 
pure state, simplex} \\
\noindent \msccode{Primary 13J30, 14P99; secondary 15B48, 52A05}

\begin{abstract}
P\'olya's Positivstellensatz and Handelman's Positivstellensatz
	are known to be concrete instances of the abstract Archimedean Representation Theorem 
		for (commutative unital) rings.
We generalise the Archimedean Representation Theorem to 
	modules over rings.
For example, consider the module of all symmetric matrices with entries in a polynomial ring,
	also known as \emph{matrix polynomials}.
P\'olya's Positivstellensatz and Handelman's Positivstellensatz
	had been generalised by Scherer and Hol, and \Le\ and \Du\ 
		 respectively to matrix polynomials,	
	using the method of effective estimates from analysis.
We show that these two Positivstellens\"atze
	for matrix polynomials
		are concrete instances 
			of our Archimedean Representation Theorem 
				in the case of the module of symmetric matrix polynomials
					over the polynomial ring.
\end{abstract}



\section*{Introduction}

Since Berr and W\"ormann \cite{BerrWoermann2001}, 
			and even earlier from W\"ormann's thesis \cite{Woermann1998},
	it is known that several Positivstellens\"atze 
		follow from the Archimedean Representation Theorem of real algebra.
In real algebraic geometry,
		a \emph{Positivstellensatz} is the sufficiency of 
			the positivity of a polynomial on a space, usually compact, 
			for it to be representable in terms of a certificate.
A \emph{certificate} is an algebraic expression 
	that immediately witnesses the strict positivity of the polynomial on that space.		
The Positivstellens\"atze that follow from the Archimedean Representation Theorem
	are characterised by their certificates forming a so-called \lq\lq\emph{archimedean}\rq\rq\!
		subsemiring of the polynomial ring.
Given a ring $A$ (always commutative with multiplicative unit),
	a subset $S \subset A$ is a \emph{subsemiring} if it contains $0$ and $1$ and is closed under addition and multiplication.
We say that $S \subset A$ is \emph{archimedean} if $S + \Z = A$, 
	where each integer $n \in \Z$ is regarded as $n \cdot 1 \in A$. 			
These archimedean Positivstellens\"atze include those of
	P\'olya \cite{Polya1928} (reproduced in  \cite[pp.\ 57--60]{HLP1952}) and
	Handelman \cite{Handelman1988}.

The Archimedean Representation Theorem
	is a criterion for an element of a ring $A$ to lie in a module $M \subset A$ over an archimedean subsemiring $S$ of $A$.
Here, by an \emph{$S$-module}, we mean a subset $M \subset A$ 
		that contains $0$, 
		is closed under addition
		and satisfies $S M \subset M$.
The fundamental Archimedean Representation Theorem was proven and rediscovered in various versions  
		by Stone, Krivine, Kadison and Dubois, among others.
Krivine's version is definitive \cite{Krivine64a, Krivine64b}.		
Prestel and Delzell gave an account of its history \cite[Section 5.6]{PrestelDelzell2001}.

When $A = \R[x_1, \dots, x_\dimension]$ is the real polynomial algebra,
	the Archimedean Representation Theorem specialises to P\'olya's Positivstellensatz and Handelman's Positivstellensatz
		for appropriate choices of an archimedean subsemiring $S$ respectively. 
The abstract criterion in the Archimedean Representation Theorem for a polynomial $f \in A$
	to be representable as a certificate
		then reduces to the strict positivity of $f$ on the relevant compact subset of real euclidean space $\R^\dimension$.

\section*{Main results and their proofs}
	
The purpose of this note, then,
	is to generalise the Archimedean Representation Theorem 
		to a criterion for an element of a module $G$ over a ring $A$
			to lie in a subsemimodule $M \subset G$ over an archimedean subsemiring $S \subset A$.
Here, we ask the reader to take note, 
	the term \lq\lq\emph{$A$-module}\rq\rq\ is in the usual sense of an additive group
		equipped with an $A$-action $(f, s) \mapsto f \cdot s: A \times G \to G$
			satisfying the usual axioms $f \cdot (s + t) = f \cdot s + f \cdot t$,
				$(f + g) \cdot s = f \cdot s + g \cdot s$, 
				$(f g) \cdot s = f \cdot (g \cdot s)$,
				and $1 \cdot s = s$, for all $f, g \in A$, $s, t \in G$.
Then $M \subset G$ 
	is a \emph{$S$-subsemimodule} 
		if it contains $0$, 
		is closed under addition 
		and satisfies $S \cdot M \subset M$.
So a $S$-module, in the above sense as used in real algebra,
	is a $S$-subsemimodule of $A$ in our terminology, 
		where $A$ is regarded as a module over itself.

We state our main results.
Given an abelian group $G$, written additively,
	and a submonoid $M \subset G$,
		an element $u \in G$ is said to be an \emph{order unit}
			of $(G, M)$ if $M + \Z u = G$. 	
An equivalent definition of an order unit
	can be given in terms of a quasiordering $\le_M$ on $G$ (i.e.\ a reflexive and transitive binary relation on $G$) associated to $M$,
		defined by $s \le_M t$ if and only if $t - s \in M$ (for $s, t \in G$).			 	
Then $u \in M$ is an order unit of $(G ,M)$
	if, for each $s \in G$, there is a positive integer $n$ such that $s \le_M nu$.			
Order units 
	were originally named by Goodearl and Handelman 
		as \lq\lq\emph{strong units}\rq\rq\!
	in the case where $M \intersect (-M) = 0$,
		i.e. where $\le_M$ is antisymmetric, and hence a partial ordering \cite{GoodearlHandelman1976}.
By 1980, they had settled on the name \lq\lq order unit\rq\rq\ \cite{GoodearlHandelman1980}, \cite{HHL1980}.				
For example, a semiring $S$ of a ring $A$ is archimedean
	if and only if $1$ is an order unit of $(A, S)$.	

Now suppose that $G$ is a module over a ring $A$, let $M \subset G$, and let $u \in M$.
Define $\nonnegativeLocus_A(G, M, u)$ 
	as the set of all group homomorphisms $\Phi : G \to \R$ to the additive group of real numbers 
	such that $\Phi|_M \ge 0$, $\Phi(u) = 1$ and
\begin{align} \label{eq: multiplicativeLaw}
\Phi(f \cdot s) = \Phi(f \cdot u) \Phi(s) & & (\forall f\in A, s \in G).
\end{align}
Given $s \in G$, we write $s > 0$ (resp.\ $s \ge 0$) on $\nonnegativeLocus_A(G, M, u)$ 
	if $\Phi(s) > 0$ (resp.\ $\Phi(s) \ge 0$) for all $\Phi \in \nonnegativeLocus_A(G, M, u)$.

\begin{theorem}
	\label{thm: ArchimedeanRepresentationTheoremForModulesOverARing}
Let $G$ be a module over a ring $A$,
	let $M \subset G$ be a subsemimodule over an archimedean subsemiring $S$ of $A$,
		and let $u \in M$ be an order unit of $(G, M)$.
Then, for each $s \in G$ 
	with $s > 0$ on $\nonnegativeLocus_A(G, M, u)$, 
		there is some positive integer $n$
			such that $n s \in M$.
\end{theorem}
	
The property that $n s \in M$ for some positive integer $n $
	witnesses that $s \ge 0$ on $\nonnegativeLocus_A(G, M, u)$
	since then $0 \le \Phi(n s) = n \Phi(s)$ 
		would imply that $\Phi(s) \ge 0$, 
			for all $\Phi \in \nonnegativeLocus_A(G, M, u)$.	
The special case of $(G, u) = (A, 1)$ is the standard Archimedean Representation Theorem in real algebra.
	This theorem was first proven algebraically by Becker and Schwartz \cite{BeckerSchwartz1983}
	(refer also to \cite[Theorem 1.5.9]{Scheiderer2009} and \cite[Theorem 6.1]{BSS12}).

However, it is desirable that the conclusion of Theorem \ref{thm: ArchimedeanRepresentationTheoremForModulesOverARing}
	be strengthened to witness the strict positivity of $s$ on $\nonnegativeLocus_A(G, M, u)$.
Let $\Q$ denote the field of rational numbers.
Let $\R$ denote the field of real numbers
	and let $\R_+$ denote the subsemifield of nonnegative real numbers.
Given a field $\K$ with $\Q \subset \K \subset \R$,
	let $\K_+ := \K \intersect \R_+$.	
Given a $\K$-algebra $A$ (always associative and commutative with multiplicative unit),
	a \emph{$\K_+$-subsemialgebra} 
		is a subsemiring $S$ of the underlying ring of $A$ 
			that contains $\K_+$.

\begin{theorem}
	\label{thm: ArchimedeanRepresentationTheoremForModulesOverAnAlgebra}
Let $\Q \subset \K \subset \R$ be a field.
Let $G$ be a module over a $\K$-algebra $A$,
	let $M \subset G$ be a subsemimodule over an archimedean 
		$\K_+$-subsemialgebra $S$ of $A$,
	and let $u \in M$ be an order unit of $(G, M)$.
Then every $s \in G$ 
	with $s > 0$ on $\nonnegativeLocus_A(G, M, u)$ 
		lies in $M$ and is, furthermore, an order unit of $(G, M)$.
\end{theorem}			

For example, when $(A,  S) = (\R, \R_+)$,
	so that $G$ is a vector space over $\R$
		and $M \subset G$ is a convex cone,
order units of $(G, M)$ 
	are known as \emph{algebraic interior points}.
This terminology is used in convex geometry \cite[III.1.6]{Barvinok2002}.		
The contrapositive of Theorem \ref{thm: ArchimedeanRepresentationTheoremForModulesOverAnAlgebra}
	then says that if a convex cone $M$ has nonempty algebraic interior (witnessed by $u$),
		then any point $s \in G$ that does not lie in the algebraic interior of $M$
			can be weakly separated from $M$ by a hyperplane.
I thank Tobias Fritz for this observation
	that Theorem \ref{thm: ArchimedeanRepresentationTheoremForModulesOverAnAlgebra} 
		specialises this hyperplane seperation theorem.
A more general hyperplane seperation theorem
	is stated in Barvinok's textbook \cite[III.1.7]{Barvinok2002}
		(see also Eidelheit \cite{Eidelheit1936} and Kakutani \cite{Kakutani1937}). 	

As promised, the conclusion that $s \in M$ is an order unit of $(G, M)$
	witnesses that $s > 0$ on $\nonnegativeLocus_A(G, M, u)$.
Indeed, for any order unit $s$ of $(G, M)$,
	there is some positive integer $n$ with $u \le_M n s$,
	hence $1 = \Phi(u) \le \Phi(n s) = n \Phi(s)$ by the monotonicity of $\Phi$,
		therefore $\Phi(s) \ge 1/n > 0$. 

The main contribution of this note is the following observation ---
	the argument of Burgdorf, Scheiderer and Schweighofer in \cite{BSS12}
		for Theorem \ref{thm: ArchimedeanRepresentationTheoremForModulesOverARing}
			in the case where $G$ is a $S$-module (in the sense of real algebra) 
				already suffices to prove Theorem \ref{thm: ArchimedeanRepresentationTheoremForModulesOverARing} in its full generality.
We recall their argument.
\begin{enumerate}[Step 1.]
\item \label{item: EHSCriterion} 
First, they recall a result of Effros, Handelman and Shen from convex geometry \cite[Theorem 1.4]{EffrosHandelmanShen1980}. 
\end{enumerate}

The result requires some terminology to state.
Let $G$ be an abelian group  written additively,
 	$M \subset G$ a submonoid, and let $u \in M$ be an order unit of $(G, M)$.
A \emph{state} of $(G, M, u)$,
	is a group homomorphism $\Phi : G \to \R$ to the additive group of real numbers
		such that $\Phi|_M \ge 0$ and $\Phi(u) = 1$.
Regard the set of states, denoted by $\stateSpace(G, M, u)$,
	as a subset of $\R^G := \prod_G \R$ via the injection $\Phi \mapsto (\Phi(g))_{g \in G} : \stateSpace(G, M, u) \hookrightarrow \R^G$.
Then $\stateSpace(G, M, u)$ is convex,
	so that we may define a \emph{pure} state 
		as an extremal point $\Phi \in \stateSpace(G, M, u)$.	
Explicitly, $\Phi \in \stateSpace(G, M, u)$ is pure if, whenever $2\Phi = \Phi_1 + \Phi_2$ for any two $\Phi_1, \Phi_2 \in \stateSpace(G, M, u)$,
	then $\Phi = \Phi_1 = \Phi_2$.

\begin{lemma} \label{lem: EHSCriterionNonStrict}
 Let $G$ be an abelian group, 
 	$M \subset G$ a submonoid, and let $u \in M$ be an order unit of $(G, M)$.
Then, for every $s \in G$ with $\Phi(s) > 0$ for all pure states $\Phi : G \to \R$ 
	of $(G ,M, u)$,
		there is some positive integer $n$ such that $n s \in M$.
\end{lemma}

In their original version of Lemma \ref{lem: EHSCriterionNonStrict},
	Effros, Handelman and Shen assumed that the quasiorder $\le_M$ is anti-symmetric,
    		or equivalently that $M \intersect (-M) = 0$.
Burgdorf, Scheiderer and Schweighofer observed that this assumption is not necessary.
They outlined a proof of Lemma \ref{lem: EHSCriterionNonStrict} 
	using two theorems from convex geometry,
			namely the Krein-Milman Theorem
				(see \cite[III.4.1]{Barvinok2002})
			and the previously mentioned hyperplane seperation theorem. 
	
\begin{enumerate}[Step 1.]
\addtocounter{enumi}{1} 
\item \label{item: formOfPureStatesOfRings} 
	Burgdorf, Scheiderer and Schweighofer showed that 
		every pure state of $(A, S, 1)$ is multiplicative, or equivalently,	
			lies in $\nonnegativeLocus_A(A, S, 1)$ \cite[Corollary 4.4]{BSS12}.
\item \label{item: substitution} 
	Therefore, the special case of Theorem \ref{thm: ArchimedeanRepresentationTheoremForModulesOverARing} when $(G, u) = (A, 1)$
		follows by applying Step \ref{item: formOfPureStatesOfRings} to 
		Lemma \ref{lem: EHSCriterionNonStrict}.
\end{enumerate}		

We observe that their verbatim argument 
	proves the following generalisation of Step \ref{item: formOfPureStatesOfRings}
		to modules over a ring.
\begin{proposition}
	\label{prop: multiplicativeAlongFibres}
Let $G$ be a module over a ring $A$,
	let $M \subset G$ be a subsemimodule over an archimedean subsemiring $S$ of $A$,
		and let $u \in M$ be an order unit of $(G, M)$.
Then each pure state $\Phi$ of $(G, M, u)$
    satisfies \eqref{eq: multiplicativeLaw}.
\end{proposition}
Burgdorf, Scheiderer and Schweighofer stated 
	Proposition \ref{prop: multiplicativeAlongFibres}
		in the special case
			where $G \subset A$ is an ideal 
			\cite[Proposition 4.1]{BSS12}. 
We reproduce their proof in the Appendix,
	where the reader may check that their argument is valid without the assumption that $G$ is contained in $A$.
	
Thus we are ready to prove Theorem \ref{thm: ArchimedeanRepresentationTheoremForModulesOverARing}.	
\begin{proof}[Proof of Theorem \ref{thm: ArchimedeanRepresentationTheoremForModulesOverARing}]
Apply Lemma \ref{lem: EHSCriterionNonStrict} to Proposition \ref{prop: multiplicativeAlongFibres}.
\end{proof}
For Theorem \ref{thm: ArchimedeanRepresentationTheoremForModulesOverAnAlgebra}, 
	we will use the following version of Lemma \ref{lem: EHSCriterionNonStrict},
		which is a special case of \cite[Corollary 2.7]{BSS12}.
\begin{lemma} \label{lem: EHSCriterionStrict}
Let $\Q \subset \K \subset \R$ be a field.
Let $G$ be a $\K$-vector space, let $M \subset G$ a $\K_+$-subsemimodule, 
 	and let $u \in M$ be an order unit of $(G, M)$.
Then every $s \in G$ with $\Phi(s) > 0$ for all pure $\Phi \in \stateSpace(G, M, u)$
	lies in $M$ and is, furthermore, an order unit of $(G, M)$.		
\end{lemma} 
\begin{proof}[Proof of Theorem \ref{thm: ArchimedeanRepresentationTheoremForModulesOverAnAlgebra}]
Apply Lemma \ref{lem: EHSCriterionStrict} to Proposition \ref{prop: multiplicativeAlongFibres}.
\end{proof}

%

\section*{Applications to matrix polynomials}
We end with some initial consequences of Theorem \ref{thm: ArchimedeanRepresentationTheoremForModulesOverAnAlgebra}.
Recall from the introduction that 
	both the Positivstellens\"atze of Handelman \cite{Handelman1988} 
												and P\'olya \cite{Polya1928} for polynomials
	are concrete instances of the regular Archimedean Representation Theorem.
These two Positivstellens\"atze were generalised by	
\Le\ and \Du\ \cite[Theorem 3]{LeDu2018}, and Scherer and Hol \cite[Theorem 3]{SchererHol2006} 
	respectively 
		to symmetric matrices whose entries are polynomials,
			also known as \emph{symmetric matrix polynomials}.
We show that these Positivstellens\"atze for matrix polynomials
		are instances of Theorem \ref{thm: ArchimedeanRepresentationTheoremForModulesOverAnAlgebra}.	

A \emph{multi-index} is a $\numberOfConstraints$-tuple
 $k = (k_1, \dots, k_\numberOfConstraints)$ of nonnegative integers. 
We use the notation $|k| := k_1 + \cdots + k_\numberOfConstraints$
		and $\ell^k := {\ell_1}^{k_1} \cdots {\ell_\numberOfConstraints}^{k_\numberOfConstraints}$,
			for any $\numberOfConstraints$-tuple 
				of polynomials $(\ell_1, \dots, \ell_\numberOfConstraints)$
					in $\R[x] := \R[x_1, \dots, x_d]$.
For any $f \in \R[x]$ and any matrix $\matrix{P} = (p_{ij})$ with real entries,
	the entrywise product $f \cdot \matrix{P} = (f p_{ij})$
		is a matrix polynomial.
Given a symmetric real matrix $\matrix{A}$,
	let $\matrix{A} \succ 0$ (resp.\ $\matrix{A} \succeq 0$)
		denote that $\matrix{A}$ is positive definite (resp.\ positive semidefinite);
			whenever this notation is used, $\matrix{A}$ is understood have real entries.

\begin{example}[\Le\ and \Du] \label{eg: MatrixHandelman}
Let $P \subset \R^d$ be a polytope 
	(i.e.\ the convex hull of a finite set of points) 
	assumed to be full-dimensional 
		(i.e. the affine span of $P$ is the entire $\R^d$).
Fix a presentation of this bounded set $P$ as the intersection of halfspaces,
	say 
\begin{equation} \label{eq: HPolytopePresentation}
	P = \{ x \in \R^d :\, \ell_1(x), \dots, \ell_\numberOfConstraints(x) \ge 0\},
\end{equation}
		where $\ell_1, \dots, \ell_\numberOfConstraints$
			are polynomials in $\R[x] := \R[x_1, \dots, x_d]$ of degree $1$.

Then any symmetric matrix $\matrix{M} = (f_{ij})$ 
			having entries $f_{ij}$ in $\R[x]$
	with $\matrix{M}(x) = (f_{ij}(x)) \succ 0$ 
		for all $x \in P$ can be written as
\begin{equation} \label{eq: HandelmanMatrixCertificate}
\matrix{M} = \sum_{|k| = \kappa} \ell^k \cdot \matrix{P_k} 
\end{equation}	
for some integer $\kappa \ge 0$ 
	and some family $\{\matrix{P_k}\}_{|k| = \kappa}$ 
			of positive definite symmetric real matrices.
\end{example}

\begin{proof}
Let $A = \R[x]$, let $S = \{\sum_k a_k \ell^k \in \R[x] : \, \forall k, \, a_k \in \R_+\} \subset A$
		be the $\R_+$-subsemialgebra generated by 
				$\ell_1, \dots, \ell_\numberOfConstraints$,
	let $G = \Sym_\matrixSize(A)$ 
		be the $A$-module 
		of symmetric $\matrixSize \times \matrixSize$ matrices with entries in $A$,
	let $M = \{ \sum_i g_i \cdot \matrix{S_i} : \, \forall i, \, g_i \in S, \matrix{S_i} \succeq 0\}$
		be the $S$-subsemimodule 
			generated by the positive semidefinite $\matrixSize \times \matrixSize$ symmetric real matrices, 
				 	and let $u = \matrix{I_\matrixSize}$ be the $\matrixSize \times \matrixSize$ identity matrix.

We shall apply Theorem \ref{thm: ArchimedeanRepresentationTheoremForModulesOverAnAlgebra}.
In order to show that $S \subset A$ is archimedean,
	it suffices to find $c_j \in \R$ such that $c_j \pm x_j \in S$ (for $j = 1, \dots, d$),
		because $x_1, \dots, x_d$ generate $A$ \cite[Lemma 1]{BerrWoermann2001}.
 But every polynomial of degree at most $1$ in $\R[x] = A$
	that is strictly positive on $P$
		can be written as a strictly positive linear combination 
			of $\ell_1, \dots, \ell_\numberOfConstraints$
			(see Handelman's paper for a proof
			\cite[Proposition I.1(b)]{Handelman1988}).
Thus $S$ contains all polynomials of degree $1$ that are strictly positive on $P$.
In particular, for all $j = 1, \dots, d$,
	the degree $1$ polynomial $c_j \pm x_j$ lies in $S$ for sufficiently large $c_j \in \R_+$
	by the compactness of $P$.
Therefore $S$ is indeed archimdean in $A$.			
		
Then $u = \matrix{I_\matrixSize}$ is an order unit of $(G, M)$ due to the following identity:
\begin{equation}
mn \matrix{I_\matrixSize} - f \cdot \matrix{A}
= \frac{1}{2}\big(
 (m  - f) \cdot (\constantBound  \matrix{I_\matrixSize} + \matrix{A}) 
 + (m  + f) \cdot (\constantBound \matrix{I_\matrixSize} - \matrix{A})
\big),
\end{equation}
which holds for all positive integers $m, n$, all $f \in A$, and all $\matrixSize \times \matrixSize$ symmetric real matrices $\matrix{A}$.
Here, note that 	the $f \cdot \matrix{A}$'s, 
				where $f \in A$ and $\matrix{A}$ ranges over the $\matrixSize \times \matrixSize$ symmetric real matrices,
	generate $G$ linearly over $\R$.
This argument is standard, see \cite[Proof of Lemma 1]{BerrWoermann2001}.	

We proceed to characterise $\nonnegativeLocus_A(G, M, u)$.
Let an arbitrary group homomorphism $\Phi : G \to \R$ in $\nonnegativeLocus_A(G, M, u)$ be given.
The $\Q$-linearity of $\Phi$ follows from the additivity, 
	by a standard argument.
In fact, $\Phi$ is $\R$-linear.
To see this, \eqref{eq: multiplicativeLaw}
	gives $\Phi(cs) = \Phi(cu)\Phi(s)$ for any $c \in \R$, $s \in G$,
		so it suffices to show that $\Phi(cu) = c\Phi(u)$.
But given any two $q_1, q_2 \in \Q$ with $q_1 \le c \le q_2$,
	the inequality chain
		$q_1 u \le_M c u \le_M q_2 u$
			implies, by the monotonicity of $\Phi$ with respect to $\le_M$, that
\begin{equation}
q_1 \Phi(u) = \Phi(q_1 u) 
	\le \Phi(cu) \le \Phi(q_2 u) = q_2\Phi(u),
\end{equation}						
thus $\Phi(cu) = c\Phi(u)$ by letting $q_1$ approach $c$ from the left
and letting $q_2$ approach $c$ from the right.

The induced map $f \mapsto \Phi(f \cdot \matrix{I_\matrixSize}) : A \to \R$ is a $\R$-algebra homomorphism
	since \eqref{eq: multiplicativeLaw} amounts to its multiplicativity
		and $\Phi(1 \cdot \matrix{I_\matrixSize}) = \Phi(\matrix{I_\matrixSize}) = 1$.	
Thus $\Phi(f \cdot \matrix{I_\matrixSize}) = f(x)$ for some $x \in \R^d$ (for $f \in A$).
In fact $x \in P$
	because $\ell_i(x) = \Phi(\ell_i \cdot \matrix{I_\matrixSize}) \ge 0$
		since $\ell_i \cdot \matrix{I_\matrixSize} \in M$,
			for all $i = 1, \dots, \numberOfConstraints$.
Now denote by $\Sym_\matrixSize(\R)$ the $\R$-vector space of $\matrixSize \times \matrixSize$ symmetric real matrices.
The restriction $\Phi|_{\Sym_\matrixSize(\R)} : \Sym_\matrixSize(\R) \to \R$
	is nonnegative-valued on the convex cone of positive semidefinite real matrices 
			(because this cone is contained in $M$)
			and hence is represented by some positive semidefinite real matrix $\matrix{S}$.
That is to say, $\Phi(\matrix{A}) = \langle \matrix{A}, \matrix{S} \rangle = \tr(\matrix{A} \matrix{S})$ 
						for $\matrix{A} \in \Sym_\matrixSize(\R)$,
					where $\langle \ ,\ \rangle$ is the Hadamard inner product
						given in terms of the trace.	
Note that $\tr(\matrix{S}) = \tr(\matrix{I_\matrixSize} \matrix{S}) = \Phi(\matrix{I_\matrixSize}) = 1$.						
Therefore, using \eqref{eq: multiplicativeLaw} again,
	for any $f \in A$ and $\matrix{A} \in \Sym_\matrixSize(\R)$,
\begin{align}
\Phi(f \cdot \matrix{A}) = \Phi(f \cdot \matrix{I_\matrixSize}) \Phi(\matrix{A})
 = f(x) \tr(\matrix{A} \matrix{S}) = \tr(f(x) \matrix{A} \matrix{S})
 = \tr( (f \cdot \matrix{A})(x) \matrix{S}).
\end{align}							
Finally, since the $f \cdot \matrix{A}$'s generate $G$ linearly over $\R$, as $f$ ranges over $A$ and $\matrix{A}$ ranges over $\Sym_\matrixSize(\R)$,
	therefore							 		
\begin{align} \label{eq: expressionForPhi}
\Phi(\matrix{M}) = \tr (\matrix{M}(x) \matrix{S}) & & (\forall \matrix{M} \in G).
\end{align}

To complete the proof, let $\matrix{M} \in G$ be given
	with $\matrix{M}(x) \succ 0$ 
		for all $x \in P$.
For any $\Phi \in \nonnegativeLocus_A(G, M, u)$,
	let $x \in P$ be the associated point 
		and $\matrix{S}$ be the associated positive semidefinite matrix with $\tr(\matrix{S}) = 1$
			such that $\Phi$ is given by \eqref{eq: expressionForPhi}.		
Then $\matrix{M}(x)^{1/2} \matrix{S} \matrix{M}(x)^{1/2}$
				is a nonzero positive semidefinite symmetric matrix,
			thus $\Phi(\matrix{M}) = \tr (\matrix{M}(x) \matrix{S}) = \tr(\matrix{M}(x)^{1/2} \matrix{S} \matrix{M}(x)^{1/2}  ) > 0$.
Since $\Phi \in \nonnegativeLocus_A(G, M, u)$ is arbitrary,
	therefore $\matrix{M}$ is an order unit of $(G ,M)$ by Theorem \ref{thm: ArchimedeanRepresentationTheoremForModulesOverAnAlgebra}.
Hence $\matrix{I_\matrixSize} \le_M n \matrix{M}$
	for some positive integer $n$, 
so that there are $a_{i, j} \in \R_+$, $\matrix{S_i} \succeq 0$	
	such that
$\matrix{M} = (1/n) (\matrix{I_\matrixSize} + \sum_i \sum_{j} a_{i,j} \ell^j \cdot \matrix{S_i})$.
Now, express the multiplicative unit $1 \in A$, which is a degree-$0$ polynomial strictly positive on $P$,
	as a strictly positive linear combination of $\ell_1, \dots, \ell_\numberOfConstraints$.
Say $1 = b_1 \ell_1 + \cdots + b_\numberOfConstraints \ell_\numberOfConstraints$, 
	where $b_1, \dots, b_\numberOfConstraints > 0$.
Therefore, for all integers $\kappa \ge \max_j |j|$,
\begin{align}
\matrix{M} &= 
\frac{1}{n}\left(
(b_1\ell_1 + \cdots + b_\numberOfConstraints \ell_\numberOfConstraints)^\kappa \cdot \matrix{I_\matrixSize}
	+ \sum_i \sum_j a_{i, j} \ell^j (b_1\ell_1 + \cdots + b_\numberOfConstraints \ell_\numberOfConstraints)^{\kappa - |j|}
										 \cdot \matrix{S_i} 
										 \right) \\
&= \frac{1}{n}\sum_{|k| = \kappa} \ell^k \cdot \left(
\binom{\kappa}{k} b^k \matrix{I_\matrixSize}
+ \sum_i \matrix{S_i} \sum_{j + p = k} a_{i, j} \binom{\kappa - |j|}{p} b^p
\right),										 
\end{align}
where $\binom{\kappa}{k} = \kappa!/(k_1! \cdots k_{\numberOfConstraints}!)$
	is the multinomial coefficient.	
Therefore $\matrix{M}$ has the desired representation as in \eqref{eq: HandelmanMatrixCertificate},
	for all integers $\kappa \ge \max_j |j|$.
Explicitly, $\matrix{P_k} := (1/n)(\binom{\kappa}{k} b^k \matrix{I_\matrixSize}
+ \sum_i \matrix{S_i} \sum_{j + p = k} a_{i, j} \binom{\kappa - |j|}{p} b^p) 
\succeq (1/n)\binom{\kappa}{k} b^k \matrix{I_\matrixSize} \succ 0$ for all $|k|  = \kappa$. 
\end{proof}

The case of $1 \times 1$ matrices in Example \ref{eg: MatrixHandelman}
	is the classical Positivstellensatz of Handelman.
\Le\ and \Du\ also gave an effective upper bound of the least $\kappa$
	required in \eqref{eq: HandelmanMatrixCertificate}
		using the corresponding effective estimates of Powers and Reznick
			for Handelman's Positivstellensatz \cite{PowersReznick2001}.

Choosing $P \subset \R^d$ as a $d$-dimensional simplex,
  	and expressing $\matrix{M}$ in terms of the barycentric coordinates of $P$
		gives Scherer and Hol's generalisation of P\'olya's Positivstellensatz.
In the $1 \times 1$ case, the observation that 
	P\'olya's Positivstellensatz 
			is Handelman's Positivstellensatz 
				for a full-dimensional simplex
	 					in barycentric coordinates 
	 						is due to Powers and Reznick \cite{PowersReznick2001}.		
Let $\Delta^d \subset \R^{d + 1}$ denote the \emph{standard $d$-dimensional simplex},
	whose vertices are $(1, 0, \dots, 0)$, $(0, 1, 0, \dots, 0)$, \dots, $(0, \dots, 0, 1)$.
A point $X = (X_0, \dots, X_d)\in \R^{d + 1}$ lies in $\Delta^d$ 
	if and only if $X_0, \dots, X_d \ge 0$ and $X_0 + \cdots + X_d = 1$.				

\begin{example}[Scherer and Hol]

Fix an integer $\degree \ge 0$,
	and let $\matrix{A}$ be a symmetric matrix 
		whose entries are forms (i.e. homogeneous polynomials) in 
			$\R[X_0, \dots, X_d]$ of degree $\degree$.
If $\matrix{A}(X) \succ 0$ for all $X \in \Delta^d$,
	then there is some integer $\kappa \ge e$,
			such that for every integer $m \ge \kappa$,	
\begin{equation} \label{eq: matrixPolyaCertificate}
(X_0 + \cdots + X_d)^{m - \degree} \cdot \matrix{A} 
	= \sum_{|k| = m} {X_0}^{k_0} \cdots {X_d}^{k_d} \cdot \matrix{P_k} 
\end{equation}	
where $\matrix{P_k} \succ 0$ 
	for all $|k| = m$.

Scherer and Hol similarly gave an effective estimate for $\kappa$
		deduced from 
	that of Powers and Reznick for P\'olya's Positivstellensatz \cite{PowersReznick2001}.
\end{example}

Other applications of Theorem \ref{thm: ArchimedeanRepresentationTheoremForModulesOverAnAlgebra} to Positivstellens\"atze should be possible.

\paragraph{Futher directions.} 
Following an anonymous referee's suggestion on an earlier version of this note,
	one may ask to what extent the results by Burgdorf, Scheiderer and Schweighofer
		\cite[Sections 4 -- 7]{BSS12}
		may be generalised from ideals in a ring $A$
			to the more general setting of modules over $A$.
Such a generalisation might yield Stellens\"atze
	for positive semidefinite matrix polynomials.

\paragraph{Acknowledgements.} I am grateful for CheeWhye Chin's encouragement to pursue this line of research.
Thank you, Tim Netzer, 
	for giving me a chance to speak at a conference at Universit\"at Innsbruck
		on some ideas in this paper.
I would also like to thank Tobias Fritz \cite{Fritz2023}, Xiangyu Liu, Mihai Putinar, 
	Claus Scheiderer, Markus Schweighofer, Wing-Keung To, and the anonymous referees
		for discussions and input.	
Finally, without help from my better half,
	I would not have had the peace of mind to complete this note.		

\section*{Appendix: Proof of Proposition \ref{prop: multiplicativeAlongFibres}}

There is nothing essentially original due to the author below ---
	the following proof of Proposition \ref{prop: multiplicativeAlongFibres} 
		is that of Burgdorf, Scheiderer and Schweighofer \cite[p.123]{BSS12},
		and reproduced verbatim below only for the convenience of the reader.
		
The author's only contribution 
	is the observation that $G$ being contained in $A$ (and hence an ideal) 
	is never used in the proof. 
As mentioned in \textit{loc.\ cit.}, precedents of Proposition \ref{prop: multiplicativeAlongFibres} can be found in 
        the work of 
                    Bonsall, Lindenstrauss and Phelps
                    \cite[Theorem 10]{BLP66},
                    Krivine \cite[Theorem 15]{Krivine64b} and
                    Handelman
                \cite[Proposition 1.2]{Handelman85}.
 
Let $A, S, G, M, u$ be as in Proposition \ref{prop: multiplicativeAlongFibres}.
Given a map $\Phi : G \to \R$,
    we associate to each $f \in A$ satisfying $\Phi(f \cdot u) \neq 0$
        a map $\Phi_f : G \to \R$ given by
\begin{align}
    \Phi_f(s) := \frac{\Phi(f \cdot s)}{\Phi(f \cdot u)} & & (\forall s\in G).
\end{align}
The reader can verify that if $\Phi$ is a state of $(G, M, u)$ and $p \in S$ satisfies $\Phi(p \cdot u) > 0$,
    then $\Phi_p$ is also a state of $(G, M, u)$.
Furthermore, if $\Phi$ is a state and $p_1, p_2 \in S$ satisfy $\Phi(p_1 \cdot u), \Phi(p_2 \cdot u) > 0$,
                so that $p_1 + p_2 \in S$ and $\Phi((p_1 + p_2) \cdot u )> 0$,
            then $\Phi_{p_1 + p_2}$ is a proper convex combination 
                        of the states $\Phi_{p_1}$ and $\Phi_{p_2}$:
\begin{equation} \label{eq: convexCombination}
    \Phi(p_1 \cdot u) \Phi_{p_1} + \Phi(p_2 \cdot u) \Phi_{p_2}
        = \Phi((p_1 + p_2) \cdot u) \Phi_{p_1 + p_2}.
\end{equation}

\begin{proof}[Proof of Proposition \ref{prop: multiplicativeAlongFibres}]
Since $S \subset A$ is archimedean and $u$ is an order unit of $(G, M)$,
    it suffices to show that \eqref{eq: multiplicativeLaw}
        holds whenever $f \in S$ and $s \in M$.

Let $f \in S$ and $s \in M$ be given.        
Then $f \cdot u\in S \cdot M \subset M$.
Hence $\Phi(f \cdot u) \ge 0$.
There are two cases: either $\Phi(f \cdot u) = 0$ or $\Phi(f \cdot u) > 0$.

\paragraph{Case 1: $\Phi(f \cdot u) = 0$.}
Then $u$ being an order unit of $(G, M)$
	gives a positive integer $\constantBound$ such that $0 \le_M s \le_M \constantBound u$.
Since $M$ is closed under the $S$-action,
    $0 \le_M f \cdot s \le_M \constantBound f \cdot u$.
Now $\Phi$ is monotone with respect to $\le_M$, so
\begin{equation}
    0 \le \Phi(f \cdot s) \le \constantBound \Phi(f \cdot u) = 0,
\end{equation}
forcing $\Phi(f \cdot s) = 0$ 
    so that both sides of \eqref{eq: multiplicativeLaw} equals to zero in this case.
    
\paragraph{Case 2: $\Phi(f \cdot u) > 0$.}
Since $S \subset A$ is archimedean, 
	there is some positive integer $\constantBound$ such that $\constantBound - f \in S$.
But increasing 	$\constantBound$ if necessary,
	we may suppose that $\constantBound > \Phi(f \cdot u)$.
Then
\begin{equation}
\Phi((\constantBound - f) \cdot u) 
	= \constantBound\Phi(u) - \Phi(f \cdot u) 
	= \constantBound - \Phi(f \cdot u) > 0.
\end{equation}
Since $f, \constantBound - f\in S$ and $\Phi(f \cdot u), \Phi((\constantBound - f) \cdot u ) >0$,
    we may apply \eqref{eq: convexCombination}
        to conclude that $\Phi_\constantBound$ is a proper convex combination of $\Phi_f$ and $\Phi_{\constantBound - f}$.
But $\Phi_\constantBound = \Phi$ (by direct calculation, using the fact that $\constantBound$ is a scalar),
    so the purity of $\Phi$ implies that $\Phi_f = \Phi$,
            which is just \eqref{eq: multiplicativeLaw}.
\end{proof}


\end{document}